\ifx\documentclass\undefined
\documentstyle[12pt]{article}
\else
\documentclass[12pt]{article}
\usepackage{latexsym}
\usepackage{amsfonts}
\usepackage{amsmath}
\usepackage{amssymb}
\fi

\author{K. Bezdek \thanks{Partially supported by a Natural Sciences and 
Engineering Research Council of Canada Discovery Grant.} (Univ. of Calgary) }

\font\tenBbb=msbm10 at 12pt         \font\sevenBbb=msbm9    \font\fiveBbb=msbm7
\newfam\Bbbfam
\textfont\Bbbfam=\tenBbb \scriptfont\Bbbfam=\sevenBbb
\scriptscriptfont\Bbbfam=\fiveBbb

\def\eps{\varepsilon}

\newtheorem{theorem}{Theorem}[section]

\newtheorem{cor}[theorem]{Corollary}
\newtheorem{con}[theorem]{Conjecture}

\newtheorem{prob}[theorem]{Problem}

\title{From the Kneser-Poulsen conjecture to ball-polyhedra
\footnote{Keywords: volume, contraction, Voronoi diagram, the Kneser-Poulsen conjecture, Schl\"afli differential formula, billiards, illumination.  
2000 Mathematical Subject Classification. Primary: 52A38, 52A40
Secondary: 52A20, 52C99}}

\begin{document}

\maketitle

\date

\begin{abstract}
A very fundamental geometric problem on finite systems of spheres was independently phrased by Kneser (1955) and Poulsen (1954). According to their well-known conjecture if a finite set of balls in Euclidean space is repositioned so that the distance between the centers of every pair of balls is decreased, then the volume of the union (resp., intersection) of the balls is decreased (resp., increased). In the first half of this paper we survey the state of the art of the Kneser-Poulsen conjecture in Euclidean, spherical as well as hyperbolic spaces with the emphases being on the Euclidean case. Based on that it seems very natural and important to study the geometry of intersections of finitely many congruent balls from the viewpoint of discrete geometry in Euclidean space. We call these sets ball-polyhedra. In the second half of this paper we survey a selection of fundamental results known on ball-polyhedra. Besides the obvious survey character of this paper we want to emphasize our definite intention to raise quite a number of open problems to motivate further research. 
\end{abstract}

\section{Introduction}
\label{zero}
The Kneser-Poulsen conjecture raises an important fundamental problem on volume measure. This paper suveys the major developments regarding this problem and orients the attention of the reader towards a number of open questions in order to generate further research progress. Anyone interested in the history of the Kneser-Poulsen conjecture as well as in the many related references is refereed to the recent paper of Bezdek and Connelly~\cite{bc} and to the elegant book of Klee and Wagon \cite{kw}.

\section{The Kneser-Poulsen conjecture}
\label{one}

Let $\|\dots \|$ denote the standard Euclidean norm of the $n$-dimensional Euclidean space ${\bf E}^n$. So, if ${\bf p}_i, {\bf p}_j$ are two points in ${\bf E}^n$, then $\|{\bf p}_i- {\bf p}_j \|$ denotes the Euclidean distance between them. It will be convenient to denote the (finite) point configuration consisting of the points ${\bf p}_1, {\bf p}_2, \dots, {\bf p}_N$ in ${\bf E}^n$ by ${\bf p}=({\bf p}_1, {\bf p}_2, \dots, {\bf p}_N)$. Now, if ${\bf p}=({\bf p}_1, {\bf p}_2, \dots, {\bf p}_N)$ and ${\bf q}=({\bf q}_1, {\bf q}_2, \dots, {\bf q}_N)$ are two configurations of $N$ points in ${\bf E}^n$ such that for all $1\le i<j\le N$ the inequality $\|{\bf q}_i- {\bf q}_j \|\le \|{\bf p}_i- {\bf p}_j \|$ holds, then we say that ${\bf q}$ is a {\it contraction} of ${\bf p}$. If ${\bf q}$ is a contraction of ${\bf p}$, then there may or may not be a continuous motion ${\bf p}(t)=({\bf p}_1(t), {\bf p}_2(t), \dots, {\bf p}_N(t))$, with  ${\bf p}_i(t)\in {\bf E}^n$ for all $0\le t\le 1$ and $1\le i\le N$ such that ${\bf p}(0)={\bf p}$ and ${\bf p}(1)={\bf q}$, and $\|{\bf p}_i(t)- {\bf p}_j(t)\|$ is monotone decreasing for all $1\le i<j\le N$. When there is such a motion, we say that ${\bf q}$ is a {\it continuous contraction} of ${\bf p}$. Finally, let $B^n({\bf p}_i, r_i)$ denote the closed $n$-dimensional ball centered at ${\bf p}_i$ with radius $r_i$ in ${\bf E}^n$ and let ${\rm Vol}_n(\dots)$ represent the $n$-dimensional volume (Lebesgue measure) in ${\bf E}^n$. In 1954 Poulsen \cite{p} and in 1955 Kneser \cite{k} independently conjectured the following for the case when $r_1=\dots=r_N$:

\medskip

\begin{con} \label{elso} If ${\bf q}=({\bf q}_1, {\bf q}_2, \dots, {\bf q}_N)$ is a contraction of ${\bf p}=({\bf p}_1, {\bf p}_2, \dots, {\bf p}_N)$ in ${\bf E}^n$, then
$${\rm Vol}_n[\cup_{i=1}^{N}B^n({\bf p}_i, r_i)]\ge {\rm Vol}_n[\cup_{i=1}^{N}B^n({\bf q}_i, r_i)].$$
\end{con}

\medskip

\begin{con} \label{masodik} If ${\bf q}=({\bf q}_1, {\bf q}_2, \dots, {\bf q}_N)$ is a contraction of ${\bf p}=({\bf p}_1, {\bf p}_2, \dots, {\bf p}_N)$ in ${\bf E}^n$, then
$${\rm Vol}_n[\cap_{i=1}^{N}B^n({\bf p}_i, r_i)]\le {\rm Vol}_n[\cap_{i=1}^{N}B^n({\bf q}_i, r_i)].$$
\end{con}

Actually, M. Kneser seems to be the one who has generated a great deal of interest in the above conjectures also via private letters written to a number of mathematicians. For more details on this see for example \cite{kw}.

\section{Nearest and farthest point Voronoi diagrams}
\label{two}

For a given point configuration ${\bf p}=({\bf p}_1, {\bf p}_2, \dots, {\bf p}_N)$ in ${\bf E}^n$ and radii $r_1, r_2, \dots , r_N$ consider the following sets:

$$V_i=\{{\bf x}\in {\bf E}^n\ |\ {\rm for\  all\  {\it j},}\ \|{\bf x}-{\bf p}_i\|^2-r_i^2\le \|{\bf x}-{\bf p}_j\|^2-r_j^2\},$$

$$V^i=\{{\bf x}\in {\bf E}^n\ |\ {\rm for\  all\  {\it j},}\ \|{\bf x}-{\bf p}_i\|^2-r_i^2\ge \|{\bf x}-{\bf p}_j\|^2-r_j^2\}.$$

The set $V_i$ (resp., $V^i$) is called the {\it nearest (resp., farthest) point Voronoi cell} of the point ${\bf p}_i$. (For a detailed discussion on nearest as well as farthest point Voronoi cells we refer the interested reader to \cite{e} and \cite{s}.) We now restrict each of these sets as follows:

$$V_i(r_i)=V_i\cap B^n({\bf p}_i, r_i),$$

$$V^i(r_i)=V^i\cap B^n({\bf p}_i, r_i).$$

We call the set $V_i(r_i)$ (resp., $V^i(r_i)$) the {\it nearest (resp., farthest) point truncated Voronoi cell} of the point ${\bf p}_i$. For each $i\neq j$ let $W_{ij}=V_i\cap V_j$ and $W^{ij}=V^i\cap V^j$. The sets $W_{ij}$ and $W^{ij}$ are the {\it walls} between the nearest point and farthest point Voronoi cells. Finally, it is natural to define the relevant {\it truncated walls} as follows:

$$W_{ij}({\bf p}_i, r_i)=W_{ij}\cap B^n({\bf p}_i, r_i)=$$ 
$$W_{ij}({\bf p}_j, r_j)=W_{ij}\cap B^n({\bf p}_j, r_j),$$

$$W^{ij}({\bf p}_i, r_i)=W^{ij}\cap B^n({\bf p}_i, r_i)=$$ 
$$W^{ij}({\bf p}_j, r_j)=W^{ij}\cap B^n({\bf p}_j, r_j).$$

\section{Csik\'os's formula}
\label{three}

The following formula discovered by Csik\'os \cite{cs1} proves Conjecture~\ref{elso} as well as Conjecture~\ref{masodik} for continuous contractions in a straighforward way in any dimension. (Actually, the planar case of the Kneser-Poulsen conjecture under continuous contractions have been proved independently in \cite{b}, \cite{cs0}, \cite{c} and \cite{bs}.)

\medskip

\begin{theorem} \label{csikos} {\it Let $n\ge 2$ and let ${\bf p}(t), 0\le t\le 1$ be a smooth motion of a point configuration in ${\bf E}^n$ such that for each $t$, the points of the configuration are pairwise distinct.
Then 
$$\frac{d}{dt}{\rm Vol}_n[\cup_{i=1}^{N}B^n({\bf p}_i(t), r_i)]=$$
$$\sum_{1\le i< j\le N} (\frac{d}{dt}d_{ij}(t))\cdot{\rm Vol}_{n-1}[W_{ij}({\bf p}_i(t), r_i)],$$

$$\frac{d}{dt}{\rm Vol}_n[\cap_{i=1}^{N}B^n({\bf p}_i(t), r_i)]= $$
$$\sum_{1\le i< j\le N} -(\frac{d}{dt}d_{ij}(t))\cdot{\rm Vol}_{n-1}[W^{ij}({\bf p}_i(t), r_i)],$$

where $d_{ij}(t)=\|{\bf p}_i(t)-{\bf p}_i(t)\|$.}
\end{theorem}

\medskip

On the one hand, Csik\'os \cite{cs2} managed to generalize his formula to configurations of balls called flowers which are sets obtained from balls with the help of operations $\cap$ and $\cup$. This work extends to hyperbolic as well as spherical space. On the other hand, Csik\'os \cite{cs3} has succeded to prove a Schl\"afli-type formula for polytopes with curved faces lying in pseudo-Riemannian Einstein manifolds, which can be used to provide another proof of Conjecture~\ref{elso} as well as Conjecture~\ref{masodik} for continuous contractions (for more details see \cite{cs3}).

\section{A short outline of the proof of Bezdek and Connelly of Conjectures~\ref{elso} and \ref{masodik} in ${\bf E}^2$}
\label{four}

In the recent paper \cite{bc} Bezdek and Connelly proved Conjecture~\ref{elso} as well as Conjecture~\ref{masodik} in the Euclidean plane. In fact, the paper contains a proof of an extension of these conjectures to flowers as well. In what follows we give an outline of the three step proof published in \cite{bc} by phrasing it through a sequence of theorems each being higher dimensional. The proofs of these results are based on the underlying Voronoi diagrams.

\medskip

\begin{theorem} \label{I} Consider $N$ moving closed $n$-dimensional balls $B^n({\bf p}_i(t), r_i)$ with $1\le i\le N, 0\le t\le 1$ in ${\bf E}^n$. If $F_i(t)$ is the contribution of the $i$th ball to the boundary of the union $\cup_{i=1}^{N}B^n({\bf p}_i(t), r_i)$ (resp., of the intersection $\cap_{i=1}^{N}B^n({\bf p}_i(t), r_i)$), then 
$$\sum_{1\le i\le N} \frac{1}{r_i}\cdot {\rm Vol}_{n-1}(F_i(t))$$
decreases (resp., increases) in t under any analytic contraction ${\bf p}(t)$ of the center points, where $0\le t\le 1$.
\end{theorem}

\medskip

\begin{theorem} \label{II} Let the centers of the closed $n$-dimensional balls $B^n({\bf p}_i, r_i)$, $1\le i\le N$ lie in the $(n-2)$-dimensional (affine) subspace $L$ of ${\bf E}^n$. If $F_i$ stands for the contribution of the $i$th ball to the boundary of the union $\cup_{i=1}^{N}B^n({\bf p}_i, r_i)$ (resp., of the intersection $\cap_{i=1}^{N}B^n({\bf p}_i, r_i)$), then 
$$\frac{1}{2\pi}\sum_{1\le i\le N} \frac{1}{r_i}\cdot {\rm Vol}_{n-1}(F_i)$$
is equal to the volume of $\cup_{i=1}^{N}B^{n-2}({\bf p}_i, r_i)$ (resp., $\cap_{i=1}^{N}B^{n-2}({\bf p}_i, r_i)$) lying $L$.
\end{theorem}

\medskip

\begin{theorem} \label{III} If ${\bf q}=({\bf q}_1, {\bf q}_2, \dots, {\bf q}_N)$ is a contraction of ${\bf p}=({\bf p}_1, {\bf p}_2, \dots, {\bf p}_N)$ in ${\bf E}^n$, then there is an analytic contraction of ${\bf p}$ onto ${\bf q}$ in ${\bf E}^{2n}$.
\end{theorem}

\medskip

Note that Theorem~\ref{I}, \ref{II} and \ref{III} imply in a straighforward way that Conjecture~\ref{elso} as well as Conjecture~\ref{masodik} hold in the Euclidean plane. Also, it is worth mentioning that somewhat surprisingly Theorem~\ref{III} (also called the leapfrog lemma) cannot be improved namely, it has been proved in \cite{mbc} that there exist point configurations
${\bf q}$ and ${\bf p}$ in ${\bf E}^n$, constructed actually in the way as it was suggested in \cite{bc}, such that ${\bf q}$ is a contraction of ${\bf p}$ in ${\bf E}^n$ and there is no continuous contraction from  ${\bf p}$ to ${\bf q}$ in ${\bf E}^{2n-1}$.

\section{Further results obtained from the proof of Bezdek and Connelly}
\label{five}

It is worth listing two additional results obtained from the proof published in \cite{bc} in order to describe a more complete picture of the status of the Kneser-Poulsen conjecture. For more details see \cite{bc}.

\medskip

\begin{theorem} Let ${\bf p}=({\bf p}_1, {\bf p}_2, \dots, {\bf p}_N)$ and ${\bf q}=({\bf q}_1, {\bf q}_2, \dots, {\bf q}_N)$ be two point configurations in ${\bf E}^n$ such that ${\bf q}$ is a piecewise-analytic contraction of ${\bf p}$ in ${\bf E}^{n+2}$. Then the conclusions of Conjecture~\ref{elso} as well as Conjecture~\ref{masodik} hold in ${\bf E}^n$.
\end{theorem}

The following generalizes a result of Gromov in \cite{g}, who proved it in the case $N\le n+1$.

\begin{theorem} If ${\bf q}=({\bf q}_1, {\bf q}_2, \dots, {\bf q}_N)$ is an arbitrary contraction of 
${\bf p}=({\bf p}_1, {\bf p}_2, \dots, {\bf p}_N)$ in ${\bf E}^n$ and $N\le n+3$, then both Conjecture~\ref{elso} and Conjecture~\ref{masodik} hold.
\end{theorem}

As a next step it would be natural to investigate the case $N=n+4$.

\medskip

\section{Kneser-Poulsen-type results for spherical and hyperbolic convex polytopes}
\label{six}

It is somewhat surprising that in spherical space for specific radius of balls (i.e. spherical caps) one can find a proof of both Conjecture~\ref{elso} and Conjecture~\ref{masodik} in all dimensions. The magic radius is $\frac{\pi}{2}$ and the following theorem describes the desired result in details. 

\medskip

\begin{theorem} \label{gombi} If a finite set of closed $n$-dimensional balls of radius $\frac{\pi}{2}$ (i.e. of closed hemispheres) in the $n$-dimensional spherical space is rearranged so that the (spherical) distance between each pair of centers does not increase, then the (spherical) $n$-dimensional volume of the intersection does not decrease and the (spherical) $n$-dimensional volume of the union does not increase. 
\end{theorem}

\medskip
The method of the proof published by Bezdek and Connelly in \cite{bc04} can be described as follows. First, one can use a leapfrog lemma to move one configuration to the other in an analytic and monotone way, but only in higher dimensions. Then the higher-dimensional balls have their combined volume (their intersections or unions) change monotonically, a fact that one can prove using Schl\"afli's differential formula. Then one can apply an integral formula to relate the volume of the higher dimensional object to the volume of the lower-dimensional object, obtaining the volume inequality for the more general discrete motions.

The following statement is a corollary of Theorem~\ref{gombi} (for details see \cite{bc04}) the Euclidean part of which has been proved independently by Alexander \cite{a85}, Capoyleas and Pach \cite{cp} and Sudakov \cite{su}.

\begin{theorem}
Let ${\bf p}=({\bf p}_1, {\bf p}_2, \dots, {\bf p}_N)$ be $N$ points on a hemisphere of the $2$-dimensional spherical space ${\bf S}^2$ (resp., points in ${\bf E}^2$), and let ${\bf q}=({\bf q}_1, {\bf q}_2, \dots, {\bf q}_N)$ be a contraction of ${\bf p}$ in ${\bf S}^2$ (resp., in ${\bf E}^2$). Then the perimeter of the convex hull of ${\bf q}$ is less than or equal to the perimeter of the convex hull of ${\bf p}$.
\end{theorem}

We remark that Theorem~\ref{gombi} extends to flowers as well moreover, a positive answer to the following problem would imply that both Conjecture~\ref{elso} and Conjecture~\ref{masodik} hold for circles in ${\bf S}^2$ (for more details on this see \cite{bc04}).

\begin{prob}
Suppose that ${\bf p}=({\bf p}_1, {\bf p}_2, \dots, {\bf p}_N)$ and ${\bf q}=({\bf q}_1, {\bf q}_2, \dots, {\bf q}_N)$ are two configurations in ${\bf S}^2$. Then prove or disprove that there is a monotone piecewise-analytic motion from
${\bf p}=({\bf p}_1, {\bf p}_2, \dots, {\bf p}_N)$ to ${\bf q}=({\bf q}_1, {\bf q}_2, \dots, {\bf q}_N)$ in ${\bf S}^4$.
\end{prob}

Note that in fact, Theorem~\ref{gombi} states a volume inequality between two spherically convex polytopes satisfying some metric conditions. The following problem searches for a natural analogue of that in hyperbolic $3$-space. In order to state it properly we recall the following. Let $A$ and $B$ be two planes in the hyperbolic $3$-space and let $A^+$ (resp., $B^+$) denote one of the two closed halfspaces bounded by $A$ (resp., $B$) such that the set $A^+\cap B^+$ is nonempty. Recall that either $A$ and $B$ intersect or $A$ is parallel to $B$ or $A$ and $B$ have a line perpendicular to both of them. Now, "the dihedral angle $A^+\cap B^+$" means not only the set in question but, also it refers to the standard angular measure of the corresponding angle between $A$ and $B$ in the first case, it refers to $0$ in the second case, and finally, in the third case it refers to the negative of the distance between $A$ and $B$ as well.

\begin{prob}
Let $P$ and $Q$ be compact convex polyhedra of the $3$-dimensional hyperbolic space with $P$ (resp., $Q$) being the intersection of the closed halfspaces $H_1^P, H_2^P, \dots, H_N^P$ (resp., $H_1^Q, H_2^Q, \dots, H_N^Q$). Assume that the dihedral angle $H_i^Q\cap H_j^Q$ is at least as large as the corresponding dihedral angle $H_i^P\cap H_j^P$ for all $1\le i<j\le N$. Then prove or disprove that the volume of $P$ is at least as large as the volume of $Q$.
\end{prob}

Using Andreev's version \cite{an} of the Koebe-Andreev-Thurston theorem and Schl\"afli's differential formula Bezdek \cite{b05} proved the following partial analogue of Theorem~\ref{gombi} in hyperbolic $3$-space.

\begin{theorem}
Let $P$ and $Q$ be nonobtuse-angled compact convex polyhedra of the same simple combinatorial type in hyperbolic $3$-space. If each inner dihedral angle of $Q$ is at least as large as the corresponding inner dihedral angle of $P$, then the volume of $P$ is at least as large as the volume of $Q$.
\end{theorem}

\medskip
\section{Alexander's conjecture}
\label{seven}
\medskip

It seems that in the Euclidean plane, for the case of the intersection of congruent disks, one can sharpen the results proved by Bezdek and Connelly \cite{bc}. Namely, Alexander \cite{a85} conjectures the following.

\begin{con} \label{alexander}
Under arbitrary contraction of the center points of finitely many congruent disks in the Euclidean plane, the perimeter of the intersection of the disks cannot decrease.
\end{con}

The analogous question for the union of congruent disks has a negative answer, as was observed by Habicht and Kneser long ago (for details see \cite{bc}). In \cite{bcc} some supporting evidence for the above conjecture of Alexander has been collected in particular, the following theorem was proved.

\begin{theorem}
Alexander's conjecture holds for continuous contractions of the center points and it holds up to $4$ congruent disks under arbitrary contractions of the center points.
\end{theorem}
 
We note that Alexander's conjecture does not hold for incronguent disks (even under continuous contractions of their center points) as it is shown in \cite{bcc}. Last but not least we remark that if Alexander's conjecture were true, then it would be a rare instance of an asymmetry between intersections and unions for Kneser-Poulsen type questions.

\medskip

\section{Disk-polygons and ball-polyhedra}
\label{eight}

The previous sections indicate a good deal of geometry on unions and intersections of balls that is worth for studying. In particular, when we restrict our attention to intersections of balls the underlying convexity suggests a broad spectrum of new analytic and combinatorial results. To make the setup ideal for discrete geometry from now on we will look at intersections of finitely many congruent closed $n$-dimensional balls with non-empty interior in ${\bf E}^n$. Also, it is natural to assume that removing any of the balls defining the ball-polyhedron in question yields the intersection of the remaining balls to become a larger set. If $n=2$, then we will call the sets in question {\it disk-polygons} and for $n\ge 3$ they will be called {\it ball-polyhedra}. This definition along with some basic properties of ball-polyhedra (resp., disk-polygons) were introduced by Bezdek in a sequence of talks at the University of Calgary in the fall of 2004. Based on that the paper \cite{blnp} written by Bezdek, L\'angi, Nasz\'odi and Papez systematically extended those investigations to get a better understanding of the geometry of ball-polyhedra (resp., disk-polygons) by proving quite a number of theorems, which one can regard the analogues of the classical theorems on convex polytopes. 

\medskip

\section{Finding the shortest billiard trajectories in disk-polygons}
\label{nine}

Billiards have been around for quite some time in mathematics and generated a great deal of research. (See for example the recent elegant book \cite{t} of Tabachnikov.)
For our purposes it seems natural to define billiard trajectories in the following way. This introduces a larger class of polygons for billiard trajectories than the traditional definition widely used in the literature. So, let $C$ be an arbitrary convex domain that is a compact convex set with non-empty interior in the Euclidean plane. Then we say that the closed polygonal path $P$ (possible with self-intersections) is a {\it generalized billiard trajectory} of $C$ if all the vertices of $P$ lie on the boundary of $C$ and if all the inner angle bisectors of $P$ are perpendicular to a supporting line of $C$ passing through the corresponding vertex of $P$. If $P$ has $N$ sides, then we say that $P$ is an {\it $N$-periodic generalized billiard trajectory} in $C$. Note that our definition of generalized billiard trajectories coincides with the traditional definition of billiard trajectories whenever the billiard table has no corner points. According to Birkhoff's well-known theorem if $B$ is a strictly convex billiard table with smooth boundary (that is if the boundary of $B$ is a simple, closed, smooth and strictly convex curve) in the Euclidean plane, then for every positive integer $N>1$ there exist (at least two) $N$-periodic billiard trajectories in $B$. This motivates the following theorem that has just been proved in \cite{bb}. In order to state that theorem in a possible short form it seems natural to introduce the following concept. Let $D$ be a disk-polygon in the Euclidean plane having the property that the pairwise distances between the centers of its generating disks of radii $r$ are at most $r$. In short, we say that $D$ is a {\it fat disk-polygon} with parameter $r>0$. In fact, it is easy to see that the disk-polygon $D$ with parameter $r$ is a fat disk-polygon if and only if the centers of the generating disks of $D$ belong to $D$ or putting it somewhat differently if and only if the center of any (closed) circular disk of radius $r$ containing $D$ belongs to $D$.

\begin{theorem}\label{altalanos}
Let $D$ be a fat disk-polygon in the Euclidean plane. Then any of the shortest generalized billiard trajectories in $D$ is a $2$-periodic one.  
\end{theorem}

Take a disk-polygon $D$ with generating disks of radii $r>0$. Then choose a positive $\eps$ not larger than the inradius of $D$ (which is the radius of the largest circular disk contained in $D$) and take the union of all circular disks of radius $\eps$ that lie in $D$. The set obtained in this way we call the {\it $\eps$-rounded disk-polygon} of $D$ and denote it by $D(\eps)$. The proof of the following theorem published in \cite{bb} is based on Theorem~\ref{altalanos}.

\begin{theorem}

Let $D$ be a fat disk-polygon in the Euclidean plane. Then any of the shortest (generalized) billiard trajectories in the $\eps$-rounded disk-polygon $D(\eps)$ is a $2$-periodic one for all $\eps>0$ being sufficiently small.    
\end{theorem}

Actually, we believe that the following even stronger statement holds (see also \cite{bb}).

\begin{con}
Let $D$ be a fat disk-polygon in the Euclidean plane. Then any of the shortest (generalized) billiard trajectories in the $\eps$-rounded disk-polygon $D(\eps)$ is a $2$-periodic one for all $\eps$ being at most as large as the inradius of $D$. 
\end{con}

Last but not least we mention the following result obtained as a corollary of Theorem~\ref{altalanos}. This might be of independent interest in particular, because it generalizes the result proved in \cite{bc89} that any closed curve of length at most $1$ can be covered by a translate of any convex domain of constant width $\frac{1}{2}$ in the Euclidean plane. As usual if $C$ is a convex domain of the Euclidean plane, then let ${\rm width}(C) $ denote the minimal width of $C$ (that is the smallest distance between two parallel supporting lines of $C$).

\begin{cor}
Let $D$ be a fat disk-polygon in the Euclidean plane. Then any closed curve of length at most $2\cdot{\rm width}(D)$ of the Euclidean plane can be covered by a translate of $D$.
\end{cor}

\medskip
It would be natural and important to look for higher dimensional analogues of these theorems.

\medskip

\section{Searching for an analogue of Steinitz theorem for standard ball-polyhedra in ${\bf E}^3$}
\label{ten}
\medskip
One can represent the boundary of a ball-polyhedron in ${\bf E}^3$ as the union of {\it vertices, edges} and {\it faces} defined in a rather natural way as follows. A boundary point is called a {\it vertex} if
it belongs to at least three of the closed balls defining the ball-polyhedron.
A {\it face} of the ball-polyhedron is the intersection of 
one of the generating closed balls with the boundary of the ball-polyhedron.
Finally, if the intersection of two faces is non-empty, then it is the
union of (possibly degenerate) circular arcs. The non-degenerate
arcs are called {\it edges} of the ball-polyhedron. 
Obviously, if a ball-polyhedron in ${\bf E}^3$ is generated by at least three balls, then it 
possesses vertices, edges and faces. Finally, a ball-polyhedron is called a  
{\it standard ball-polyhedron} if its vertices, edges and faces (together with the empty set and the ball-polyhedron itself) form an algebraic lattice with respect to containment. We note that not every ball-polyhedron of ${\bf E}^3$ is a standard one a fact, that is somewhat surprising and is responsible for some of the difficulties arising at studying ball-polyhedra in general (for more details see \cite{bn} as well as \cite{blnp}).

In this survey paper, a {\it graph} is always a non-oriented one and has finitely many
vertices and edges. Recall that a graph is {\it $3$-connected} if it has at least four vertices
and deleting any two vertices yields a connected graph. Also, a graph is called 
{\it simple} if it contains no loops (edges with identical end-points)
and no parallel edges (edges with the same two end-points). Finally, a graph is {\it planar} if it can be drawn in the Euclidean plane without crossing edges. Now, recall that according to the well-known theorem of Steinitz a graph is the edge-graph of some convex polyhedron in ${\bf E}^3$ if, and
only if, it is simple, planar and $3$-connected. As a partial analogue of Steinitz theorem for ball-polyhedra the following theorem is proved in \cite{blnp}.

\begin{theorem}
The edge-graph of any standard ball-polyhedron in ${\bf E}^3$ is a
simple, planar and $3$-connected graph.
\end{theorem}

Based on that it would be highly interesting to find an answer to the following question raised in \cite{blnp}.

\begin{prob}
Prove or disprove that every simple, planar and $3$-connected graph is
the edge-graph of some standard ball-polyhedron in ${\bf E}^3$.
\end{prob}

\medskip

\section{On global rigidity of ball-po\-ly\-hed\-ra in ${\bf E}^3$}
\label{eleven}
\medskip

One of the best known 
results on the geometry of convex polyhedra is Cauchy's 
rigidity theorem: If two convex polyhedra $P$ and $Q$ in ${\bf E}^3$ are combinatorially
equivalent with the corresponding facets being congruent, then also the angles
between the corresponding pairs of adjacent facets are equal and thus, $P$ 
is congruent to $Q$. For more details on Cauchy's rigidity theorem and on its extensions
we refer the interested reader to \cite{C}. In order to phrase properly the main theorem of this section we need to recall the following terminology. To each edge of a ball-polyhedron in ${\bf E}^3$ we can assign an
{\it inner dihedral angle}. Namely, take any point ${\bf p}$ in the
relative interior of the edge and take the two balls that contain 
the two faces of the ball-polyhedron meeting along that edge. Now, the
inner dihedral angle along this edge is the angle of the two half-spaces 
supporting the two balls at ${\bf p}$. 
The angle in question is obviously independent of the choice of ${\bf p}$.
Finally, at each vertex of a face of a ball-polyhedron there is a
{\it face angle} formed by the two edges meeting at the given vertex
(which is in fact, the angle between the two tangent half-lines of the two edges
meeting at the given vertex). We say that the standard ball-polyhedron $P$ in ${\bf E}^3$ is {\it globally rigid with respect
to its face angles} (resp. {\it its inner dihedral angles}) if the following holds:
If $Q$ is another standard ball-polyhedron in ${\bf E}^3$ whose face-lattice is isomorphic 
to that of $P$ and whose face angles (resp. inner dihedral angles)
are equal to the corresponding face angles (resp. inner dihedral angles)
of $P$, then $Q$ is similar to $P$. (Note that in case the family of ball-polyhedra is defined with the additional restriction that the radii of the generating balls are all equal to say, $1$, then in the above definition of global rigidity "similar" should be replaced by "congruent" as in \cite{bn}.)  A ball-polyhedron of ${\bf E}^3$ is called {\it triangulated} if all its faces are bounded by three edges. It is easy to see that any triangulated ball-polyhedron is in fact, a standard one. The following theorem has been proved in \cite{bn}.

\begin{theorem}
Let $P$ be a triangulated ball-polyhedron in ${\bf E}^3$. 
Then $P$ is globally rigid with respect to its face angles.
\end{theorem}

It remains to be a challanging problem to answer the following related question.

\begin{prob}
Let $P$ be a triangulated ball-polyhedron in ${\bf E}^3$. 
Prove or disprove that $P$ is globally rigid with respect to its dihedral angles.
\end{prob}

Finally, we mention that one can regard the above problem as an analogue of Stokker's conjecture \cite{st} according to which for convex polyhedra the face-lattice and the dihedral angles determine the face angles.

\section{Illumination of ball-polyhedra}
\label{twelve}

As we have mentioned before \cite{blnp} lays a broad ground for future study of ball-polyhedra by proving several new properties of them and raising open research problems as well. This list includes among many things analogues of the classical separation theorems of convex polytopes, a Kirchberger-type theorem, analogues of the Caratheodory theorem and the Euler-Poincare formula for ball-polyhedra. Here we want to focus on another possible direction for research. Let $K$ be a convex body (i.e. a compact convex set with nonempty interior) in the $n$-dimensional Euclidean space ${\bf E}^{n}, n\geq 2$. According to Hadwiger (see \cite{b06}) an exterior point ${\bf p}\in {\bf E}^{n}\setminus K$ of $K$ illuminates the boundary point ${\bf q}$ of $K$ if the half line emanating from ${\bf p}$ passing through ${\bf q}$ intersects the interior of $K$ (at a point not between ${\bf p}$ and ${\bf q}$). Furthermore, a family of exterior points of $K$ say, ${\bf p}_1, {\bf p}_2, \dots , {\bf p}_N$ illuminates $K$ if each boundary point of $K$ is illuminated by at least one of the point sources ${\bf p}_1, {\bf p}_2, \dots , {\bf p}_N$. Finally, the smallest $N$ for which there exist $N$ exterior points of $K$ that illuminate $K$ is called the {\it illumination number} of $K$ denoted by $I(K)$. In 1960, Hadwiger (see \cite{b06}) raised the following amazingly elementary but, very fundamental question. An equivalent but somewhat different looking concept of illumination was introduced by Boltyanski in the same year. There he proposed to use directions (i.e. unit vectors) instead of point sources for the illumination of convex bodies (for more details see \cite{b06}). Based on these circumstances the following conjecture we call the Boltyanski-Hadwiger illumination conjecture. According to this conjecture the illumination number $I(K)$ of any convex body $K$ in ${\bf E}^{n}, n\geq 2$ is at most $2^n$ and $I(K)=2^n$ if and only if $K$ is an affine $n$-cube. This conjecture is easy to prove for $n=2$ but, it is open for all $n\ge 3$. 

The following statement follows from the Separation Lemma of Bezdek \cite{b06}. In order to state it properly we need to recall two basic notions. Let $K$ be a convex body in ${\bf E}^n$ and let $F$ be a face of
$K$ that is let $F$ be the intersection of $K$ with some of its supporting hyperplanes. 
The {\it Gauss image} $\nu (F)$ of the face $F$ is the set of
all points (i.e. unit vectors) ${\bf u}$ of the $(n-1)$-dimensional unit sphere ${\bf S}^{n-1}\subset
{\bf E}^n$ centered at the origin ${\bf o}$ of ${\bf E}^n$ for which the supporting
hyperplane of $K$ with outer normal vector ${\bf u}$ contains $F.$
It is easy to see that the Gauss images of distinct faces of $K$
have disjoint relative interiors in ${\bf S}^{n-1}$ and $\nu (F)$ is compact and spherically
convex for any face $F$. Let $C\subset {\bf S}^{n-1}$ be a set of finitely many points. Then the {\it covering radius} of $C$ is the smallest positive real number $r$ with the property that the family of spherical balls of radii $r$ centered at the points of $C$ cover ${\bf S}^{n-1}$.

\begin{theorem}\label{becsles}
Let $K\subset {\bf E}^n$, $n\geq 3$ be a convex body and let $r$ be a positive real number with the property that the Gauss image $\nu (F)$ of any face $F$ of $K$ can be covered by a spherical ball of radius $r$ in ${\bf S}^{n-1}$. Moreover, assume that there exist $N$ points of ${\bf S}^{n-1}$ with covering radius $R$ satisfying the inequality $r+R\le\frac{\pi}{2}$. Then $I(K)\le N$. 
\end{theorem}

Using Theorem~\ref{becsles} as well as the optimal codes for the covering radii of four and five points on ${\bf S}^2$ (\cite{ft}) one can easily prove the result stated below. (In fact, weaker but, still reasonable estimates can be proved for larger values of $x$ (relative to $r$) by taking into account additional (optimal) codes from \cite{ft}.)

\medskip

\begin{theorem} Let $B(x, r)$ be a ball-polyhedron in ${\bf E}^{3}$ having the property that the diameter of the centers of its generating balls is at most $x$, where $0<x< 2r$ with $r$ standing for the radii of the generating balls of $B(x, r)$. Then for $0<x\le 0.57r$ we have that $I(B(x, r))=4$ and for $0.57r<x\le 0.77r$ we get that $I(B(x, r))\le 5$. 
\end{theorem}

\medskip

Based on this it is tempting to raise the following question.

\begin{prob}
Prove or disprove that if $B$ is an arbitrary ball-polyhedron of ${\bf E}^{3}$, then $I(B)\le 5$. More generally, prove or disprove that there exists a universal constant $c>0$ such that the illumination number of any $n$-dimensional ball-polyhedron in ${\bf E}^{n}$ is smaller than $(2-c)^n$ for all $n\ge 3$. 
\end{prob}

\vspace{1cm}

\medskip

\noindent
K\'aroly Bezdek,
Department of Mathematics and Statistics,
2500 University drive N.W.,
University of Calgary, AB, Canada, T2N 1N4.
\newline
{\sf e-mail: bezdek@math.ucalgary.ca}

\end{document}